\documentclass[11pt]{article}
\usepackage{amsmath}
\usepackage{amssymb}
\usepackage{latexsym}
\begin{document}
\title{Some characterizations of Hom-Leibniz algebras}
\author{ A. Nourou ISSA \\ D\'epartement de Math\'ematiques, Universit\'e
d'Abomey-Calavi, \\ 01 BP 4521 Cotonou 01, Benin. \\
E-mail: woraniss@yahoo.fr}
\date{}
\maketitle
\begin{abstract}
Some basic properties of Hom-Leibniz algebras are found. These properties are the Hom-analogue of corresponding well-known properties of Leibniz algebras. Considering the Hom-Akivis algebra associated to a given Hom-Leibniz algebra, it is observed that the Hom-Akivis identity leads to an additional property of Hom-Leibniz algebras, which in turn gives a necessary and sufficient condition for Hom-Lie admissibility of Hom-Leibniz algebras. A necessary and sufficient condition for Hom-power associativity of Hom-Leibniz algebras is also found.
\end{abstract}
\section{Introduction}
\par
The theory of Hom-algebras originated from the introduction of the notion of a Hom-Lie algebra by J.T. Hartwig, D. Larsson and S.D. Silvestrov [5] in the study of algebraic structures describing some q-deformations of the Witt and the Virasoro algebras. A Hom-Lie algebra is characterized by a Jacobi-like identity (called the Hom-Jacobi identity) which is seen as the Jacobi identity twisted by an endomorphism of a given algebra. Thus, the class of Hom-Lie algebras contains the one of Lie algebras.
\par
Generalizing the well-known construction of Lie algebras from associative algebras, the notion of a Hom-associative algebra is introduced by A. Makhlouf and S.D. Silvestrov [12] (in fact the commutator algebra of a Hom-associative algebra is a Hom-Lie algebra). The other class of Hom-algebras closely related to Hom-Lie algebras is the one of Hom-Leibniz algebras [12] (see also [8]) which are the Hom-analogue of Leibniz algebras [9]. Roughly, a Hom-type generalization of a given type of algebras is defined by a twisting of the defining identities with a linear self-map of the given algebra. For various Hom-type algebras one may refer, e.g., to [10],  [11], [15], [16], [7]. In [14] D. Yau showed a way of constructing Hom-type algebras starting from their corresponding untwisted algebras and a self-map.
\par
In [9] (see also [3], [4]) the basic properties of Leibniz algebras are given. The main purpose of this note is to point out that the Hom-analogue of some of these properties holds in Hom-Leibniz algebras (section 3). Considering the Hom-Akivis algebra associated to a given Hom-Leibniz algebra, we observe that the property in Proposition 3.3 is the expression of the Hom-Akivis identity. As a consequence we found a necessary and sufficient condition for the Hom-Lie admissibility of Hom-Leibniz algebras (Corollary 3.4). Generalizing power-associativity of rings and algebras [2], the notion of the (right) $n$th Hom-power associativity $x^{n}$ of an element $x$ in a Hom-algebra is introduced by D. Yau [17], as well as Hom-power associativity of Hom-algebras. We found that $x^{n} = 0$, $n \geq 3$, for any $x$ in a left Hom-Leibniz algebra $(L, \cdot , \alpha )$ and that $(L, \cdot , \alpha )$ is Hom-power associative if and only if $ \alpha (x) x^{2} = 0$, for all $x$ in $L$ (Theorem 3.7). Then we deduce , as a particular case, corresponding characterizations of left Leibniz algebras (Corollary 3.8). Apart of the (right) $n$th Hom-power of an element of a Hom-algebra [17], we consider in this note the left $n$th Hom-power of the given element. This allows to prove the Hom-analogue (see Theorem 3.10) of a result of D.W. Barnes ([4], Theorem 1.2 and Corollary 1.3) characterizing left Leibniz algebras. In section 2 we recall some basic notions on Hom-algebras. Modules, algebras, and linearity are meant over a ground field $\mathbb K$ of characteristic $0$.
\section{Preliminaries}
\par
In this section we recall some basic notions related to Hom-algebras. These notions are introduced in [5], [10], [12], [14], [7]. \\
\par
{\bf Definition 2.1.} A {\it Hom-algebra} is a triple $(A, \cdot , \alpha)$ in which $A$ is a $\mathbb K$-vector space, `` $\cdot$ '' a binary operation on $A$ and $\alpha : A \rightarrow A$ is a linear map (the twisting map) such that $\alpha (x \cdot y) = \alpha (x) \cdot \alpha (y)$ (multiplicativity), for all $x,y$ in $A$.\\
\par
{\bf Remark 2.2.} A more general notion of a Hom-algebra is given (see, e.g., [10], [12]) without the assumption of multiplicativity and $A$ is considered just as a $\mathbb K$-module. For convenience, here we assume that a Hom-algebra $(A, \cdot , \alpha)$ is always multiplicative and that $A$ is a $\mathbb K$-vector space. \\
\par
{\bf Definition 2.3.} Let $(A, \cdot , \alpha)$ be a Hom-algebra.
\par
(i) The {\it Hom-associator} of $(A, \cdot , \alpha)$ is the trilinear map $ as : A \times A \times A \rightarrow A$ defined by \\
\par
$as(x,y,z) = xy \cdot \alpha (z) - \alpha (x) \cdot yz$, \\
\\
for all $x,y,z$ in $A$.
\par
(ii) $(A, \cdot , \alpha)$ is said {\it Hom-associative} if $as(x,y,z) = 0$ (Hom-associativity), for all $x,y,z$ in $A$. \\
\par
{\bf Remark 2.4.} If $\alpha = Id$ (the identity map) in $(A, \cdot , \alpha)$, then its Hom-associator is just the usual associator of the algebra $(A, \cdot )$. In Definition 2.1, the Hom-associativity is not assumed, i.e.  $as(x,y,z) \neq 0$ in general. In this case $(A, \cdot , \alpha)$ is said non-hom-associative [7] (or Hom-nonassociative [14]; in [11], $(A, \cdot , \alpha)$ is also called a nonassociative Hom-algebra). This matches the generalization of associative algebras by the nonassociative ones. \\
\par
{\bf Definition 2.5.} (i) A {\it (left) Hom-Leibniz algebra} is a Hom-algebra $(A, \cdot , \alpha)$ such that the identity \\
\par
$ \alpha (x) \cdot yz = xy \cdot \alpha (z) + \alpha (y) \cdot xz$ \hfill (2.1) \\
\\
holds for all $x,y,z$ in $A$. \\
(ii) A {\it Hom-Lie algebra} is a Hom-algebra $(A, [-,-] , \alpha)$ such that the binary operation ``$[-,-]''$ is skew-symmetric and the {\it Hom-Jacobi identity} \\
\par
$J_{\alpha}(x,y,z) = 0$ \hfill (2.2) \\
\\
holds for all $x,y,z$ in $A$ and $J_{\alpha}(x,y,z) := [[x,y], \alpha (z)] + [[y,z], \alpha (x)] + [[z,x], \alpha (y)]$ is called the {\it Hom-Jacobian}. \\
\par
{\bf Remark 2.6.} The original definition of a Hom-Leibniz algebra [12] is related to the identity \\
\par
$ xy \cdot \alpha (z) = xz \cdot \alpha (y) + \alpha (x) \cdot yz$ \hfill (2.3) \\
\\
which is expressed in terms of (right) adjoint homomorphisms $Ad_{y}x := x \cdot y$ of $(A, \cdot , \alpha)$. This justifies the term of ``(right) Hom-Leibniz algebra'' that could be used for the Hom-Leibniz algebra defined in [12]. The dual of (2.3) is (2.1) and in this note we consider only left Hom-Leibniz algebras. For $\alpha = Id$ in $(A, \cdot , \alpha)$ (resp. $(A, [-,-] , \alpha)$), any Hom-Leibniz algebra (resp. Hom-Lie algebra) is a Leibniz algebra $(A, \cdot )$ [3], [9] (resp. a Lie algebra $(A, [-,-] )$). As for Leibniz algebras, if the operation ``$\cdot$'' of a given Hom-Leibniz algebra $(A, \cdot , \alpha)$ is skew-symmetric, then $(A, \cdot , \alpha)$ is a Hom-Lie algebra (see [12]). 
\par
In terms of Hom-associators, the identity (2.1) is written as \\
\par
$as(x,y,z) = - \alpha (y) \cdot xz$ \hfill (2.4) \\
\\
Therefore, from Definition 2.3 and Remark 2.4, we see that Hom-Leibniz algebras are examples of non-Hom-associative algebras. \\
\par
{\bf Definition 2.7.} [7] A {\it Hom-Akivis algebra} is a quadruple \\ $(A, [-,-] , [-,-,-], \alpha)$ in which $A$ is a vector space, ``$[-,-]$'' a skew-symmetric binary operation on $A$, ``$[-,-,-]$'' a ternary operation on $A$ and $\alpha :A \rightarrow A$ a linear map such that the {\it Hom-Akivis identity} \\
\par
$J_{\alpha}(x,y,z) = \sigma [ x,y,z ] - \sigma [ y,x,z ]$ \hfill (2.5) \\
\\
holds for all $x,y,z$ in $A$, where $\sigma$ denotes the sum over cyclic permutation of $x,y,z$. \\
\par
Note that when $\alpha = Id$ in a Hom-Akivis algebra $(A, [-,-] , [-,-,-], \alpha)$, then one gets an Akivis algebra $(A, [-,-] , [-,-,-])$. Akivis algebras were introduced in [1] (see also references therein), where they were called $W$-algebras. The term ``Akivis algebra'' for these objects is introduced in [6].
\par
In [7], it is observed that to each non-Hom-associative algebra is associated a Hom-Akivis algebra (this is the Hom-analogue of a similar relationship between nonassociative algebras and Akivis algebras [1]). In this note we use the specific properties of the Hom-Akivis algebra associated to a given Hom-Leibniz algebra to derive a property characterizing Hom-Leibniz algebras.
\section{Characterizations}
\par
In this section, Hom-versions of some well-known properties of left Leibniz algebras are displayed. Considering the specific properties of the binary and ternary operations of the Hom-Akivis algebra associated to a given Hom-Leibniz algebra, we infer a characteristic property of Hom-Leibniz algebras. This property in turn allows to give a necessary and sufficient condition for the Hom-Lie admissibility of these Hom-algebras. The Hom-power associativity of Hom-Leibniz algebras is considered.
\par
Let $(L, \cdot , \alpha)$ be a Hom-Leibniz algebra and consider on $(L, \cdot , \alpha)$ the operations \\
\par
$[x,y] := x \cdot y - y \cdot x$ \hfill (3.1) \\
\par
$ [ x,y,z ] := as(x,y,z) $. \hfill (3.2) \\
\\
Then the operations (3.1) and (3.2) define on $L$ a Hom-Akivis structure [7]. We have the following \\
\par
{\bf Proposition 3.1.} {\it Let $(L, \cdot, \alpha)$ be a Hom-Leibniz algebra. Then \\
(i) $( x \cdot y + y \cdot x) \cdot \alpha (z) = 0$, \\
(ii) $\alpha (x) \cdot [y,z] = [x \cdot y, \alpha (z)] + [\alpha (y), x \cdot z]$, \\
for all $x,y,z$ in $L$ }. \\
\\
{\bf Proof.} The identity (2.1) implies that 
\par
$ xy \cdot \alpha (z) = \alpha (x) \cdot yz - \alpha (y) \cdot xz$. \\
Likewise, interchanging $x$ and $y$, we have
\par
$ yx \cdot \alpha (z) = \alpha (y) \cdot xz - \alpha (x) \cdot yz$. \\
Then, adding memberwise these equalities above, we come to the property (i). Next we have 
\begin{eqnarray}
[x \cdot y, \alpha (z)] + [\alpha (y), x \cdot z] &=& xy \cdot \alpha (z) - \alpha (z) \cdot xy \nonumber \\
&+& \alpha (y) \cdot xz - xz \cdot \alpha (y) \; \nonumber \\
&=& \alpha (x) \cdot yz - \alpha (z) \cdot xy - xz \cdot \alpha (y) \; \mbox {(by (2.1))} \nonumber \\
&=& \alpha (x) \cdot yz - zx \cdot \alpha (y) - \alpha (x) \cdot zy \nonumber \\
&-& xz \cdot \alpha (y) \; \mbox {(by (2.1))} \nonumber \\
&=& \alpha (x) \cdot yz  - \alpha (x) \cdot zy \; \mbox {(by (i))} \nonumber \\
&=& \alpha (x) \cdot [y,z]. \nonumber
\end{eqnarray}
and so we get (ii). \hfill $\square$ \\
\par
{\bf Remark 3.2.} If set $\alpha = Id$ in Proposition 3.1, then one recovers the well-known properties of Leibniz algebras: $( x \cdot y + y \cdot x) \cdot z = 0$ and $x \cdot [y,z] = [x \cdot y,z] + [y,x \cdot z]$ (see [3], [9]). \\
\par
{\bf Proposition 3.3.} {\it Let $(L, \cdot , \alpha)$ be a Hom-Leibniz algebra. Then \\
\par
 $J_{\alpha}(x,y,z) = \sigma xy \cdot \alpha (z)$, \hfill (3.3) \\
\\
for all $x,y,z$ in $L$ }. \\
\\
{\bf Proof.} Considering (2.5) and then applying (3.2) and (2.4), we get \\ $J_{\alpha}(x,y,z) = \sigma [-\alpha (y) \cdot xz] - \sigma [-\alpha (x) \cdot yz] = \sigma [\alpha (x) \cdot yz -\alpha (y) \cdot xz] = \sigma xy \cdot \alpha (z)$ (by (2.1)). \hfill $\square$ \\
\par
One observes that (3.3) is the specific form of the Hom-Akivis identity (2.5) in case of Hom-Leibniz algebras.
\par
The skew-symmetry of the operation ``$\cdot$'' of a Hom-Leibniz algebra $(L, \cdot , \alpha)$ is a condition for $(L, \cdot , \alpha)$ to be a Hom-Lie algebra [12]. From Proposition 3.3 one gets the following necessary and sufficient condition for the Hom-Lie admissibility [12] of a given Hom-Lie algebra. \\
\par
{\bf Corollary 3.4.} {\it A Hom-Leibniz algebra  $(L, \cdot , \alpha)$ is Hom-Lie admissible if and only if $\sigma xy \cdot \alpha (z) = 0$, for all $x,y,z$ in $L$ }. \hfill $\square$ \\
\par
In [17] D. Yau introduced Hom-power associative algebras which are seen as a generalization of power-associative algebras. It is shown that some important properties of power-associative algebras are reported to Hom-power associative algebras.
\par
Let $A$ be a Hom-Leibniz algebra with a twisting linear self-map $\alpha$ and the binary operation on $A$ denoted by juxtaposition. We recall the following \\
\par
{\bf Definition 3.5.} [17] Let $x \in A$ and denote by ${\alpha}^{m}$ the $m$-fold composition of $m$ copies of $\alpha$ with ${\alpha}^{0} := Id$.
\par
(1) The $n$th Hom-power $x^{n} \in A$ of $x$ is inductively defined by \\
\par
$x^{1} = x$, \;\; $x^{n} = x^{n-1} {\alpha}^{n-2}(x)$ \hfill (3.4) \\
\\
for $n \geq 2$.
\par
(2) The Hom-algebra $A$ is {\it nth Hom-power associative} if \\
\par
$x^{n} = {\alpha}^{n-i-1}(x^{i}) {\alpha}^{i-1}(x^{n-i})$ \hfill (3.5) \\
\\
for all $x \in A$ and $i \in \{ 1,...,n-1 \}$.
\par
(3) The Hom-algebra $A$ is {\it up to nth Hom-power associative} if $A$ is $k$th Hom-power associative for all $k \in \{ 2,...,n \}$.
\par
(4) The Hom-algebra $A$ is {\it Hom-power associative} if $A$ is $n$th Hom-power associative for all $n \geq 2$. \\
\par
The following result provides a characterization of third Hom-power associativity of Hom-Leibniz algebras. \\
\par
{\bf Lemma 3.6.} {\it Let $(L, \cdot , \alpha)$ be a Hom-Leibniz algebra. Then
\par
(i) $x^{3} = 0$, for all $x \in L$; 
\par
(ii) $(L, \cdot , \alpha)$ is third Hom-power associative if and only if ${\alpha}(x)x^{2} = 0$, for all $x \in L$}. \\
\par
{\bf Proof.} From (3.4) we have $x^{3} := x^{2}{\alpha}(x)$. Therefore the assertion (i) follows from Proposition 3.1(i) if set $y=x=z$.
\par
Next, from (3.5) we note that the $i=2$ case of $n$th Hom-power associativity is automatically satisfied since this case is $x^{3} = {\alpha}^{0}(x^{2}) {\alpha}^{1}(x^{1}) = x^{2}{\alpha}(x)$, which holds by definition. The $i=2$ case says that $x^{3} = {\alpha}^{1}(x) {\alpha}^{0}(x^{2}) = {\alpha}(x)x^{2}$. Therefore, since $x^{2}{\alpha}(x) = 0$ naturally holds by Proposition 3.1 (i), we conclude that the third Hom-power associativity of $(L, \cdot , \alpha)$ holds if and only if ${\alpha}(x)x^{2} = 0$ for all $x \in L$, which proves the assertion (ii). \hfill $\square$ \\
\par
The following result shows that the condition in Lemma 3.6 is also necessary and sufficient for the Hom-power associativity of $(L, \cdot , \alpha)$. To prove this, we rely on the main result of [17] (see Corollary 5.2). \\
\par
{\bf Theorem 3.7.} {\it Let $(L, \cdot , \alpha)$ be a Hom-Leibniz algebra. Then
\par
(i) $x^{n} = 0$, $n \geq 3$, for all $x \in L$; 
\par
(ii) $(L, \cdot , \alpha)$ is Hom-power associative if and only if ${\alpha}(x)x^{2} = 0$, for all $x \in L$}. \\
\par
{\bf Proof.} The proof of (i) is by induction on $n$: the first step $n=3$ holds by Lemma 3.6(i); now if  suppose that $x^{n} = 0$, then $x^{n+1} := x^{(n+1)-1}{\alpha}^{(n+1)-2}(x) = x^{n}{\alpha}^{n-1}(x) = 0$ so we get (i).
\par
Corollary 5.2 of [17] says that, for a multiplicative Hom-algebra, the Hom-power associativity is equivalent to both of the conditions \\
\par
$x^{2}{\alpha}(x) = {\alpha}(x)x^{2}$ and $x^{4} =  {\alpha}(x^{2}){\alpha}(x^{2})$. \hfill (3.6) \\
\\
In the situation of multiplicative left Hom-Leibniz algebras, the first equality of (3.6) is satisfied by Lemma 3.6(i) and the hypothesis ${\alpha}(x)x^{2} = 0$. Next we have, from (3.5):
\par
case $i=1$: $x^{4} :=  {\alpha}^{4-2}(x){\alpha}^{0}(x^{3}) = {\alpha}^{2}(x)x^{3}$, 
\par
case $i=2$: $x^{4} := {\alpha}(x^{2}){\alpha}(x^{2})$,
\par
case $i=3$: $x^{4} := {\alpha}^{0}(x^{3}){\alpha}^{2}(x) =  x^{3}{\alpha}^{2}(x)$. \\
Because of the assertion (i) above, only the case $i=2$ is of interest here. From one side we have $x^{4} = 0$ (by (i)) and, from the other side we have ${\alpha}(x^{2}){\alpha}(x^{2}) = [{\alpha}(x)]^{2}{\alpha}(x^{2}) = 0$ (by multiplicativity and Proposition 3.1(i)). Therefore, Corollary 5.2 of [17] now applies and we conclude that (3.6) holds (i.e. $(L, \cdot , \alpha)$ is Hom-power associative) if and only if ${\alpha}(x)x^{2} = 0$, which proves (ii). \hfill $\square$ \\
\par
Let $A$ be an algebra (over a field of characteristic $0$). For an element $x \in A$, the {\it right powers} are defined by \\
\par
$x^{1} = x$, and $x^{n+1} = x^{n}x$ \hfill (3.7) \\
\\
for $n \geq 1$. Then $A$ is power-associative if and only if \\
\par
$x^{n} = x^{n-i}x^{i}$ \hfill (3.8) \\
\\
for all $x \in A$, $n \geq 2$, and $i \in \{ 1,...,n-1 \}$. By a theorem of Albert [2], $A$ is power-associative if only if it is third and fourth power-associative, which in turn is equivalent to \\
\par
$x^{2}x = xx^{2}$ and $x^{4} =  x^{2}x^{2}$. \hfill (3.9) \\
\\
for all $x \in A$.
\par
Some consequences of the results above are the following simple characterizations of (left) Leibniz algebras. \\
\par
{\bf Corollary 3.8.} {\it Let $(L, \cdot )$ be a left Leibniz algebra. Then
\par
(i) $x^{n} = 0$, $n \geq 3$, for all $x \in L$; 
\par
(ii) $(L, \cdot)$ is power-associative if and only if $xx^{2} = 0$, for all $x \in L$}. \\
\par
{\bf Proof.} The part (i) of this corollary follows from (3.7) and Theorem 3.7(i) when $\alpha = Id$ (we used here the well-known property $(xy + yx)z = 0$ of left Leibniz algebras). The assertion (ii) is a special case of Theorem 3.7(ii) (when $\alpha = Id$), if keep in mind the assertion (i), (3.8), and (3.9). \hfill $\square$ \\
\par
{\bf Remark 3.9.} Although the condition $xx^{2} = 0$ does not always hold in a left Leibniz algebra  $(L, \cdot )$, we do have $xx^{2} \cdot z = 0$ for all $x,z \in L$ (again, this follows from the property $(xy + yx)z = 0$). In fact, $b \cdot z = 0$, $z \in L$, where $b \neq 0$ is a left $m$th power of $x$ ($m \geq 2$), i.e. $b = x(x(...(xx)...))$ ([4], Theorem 1.2 and Corollary 1.3). \\
\par
Let call the {\it nth right Hom-power} of $x \in A$ the power defined by (3.4), where $A$ is a Hom-algebra. Then one may consider the {\it nth left Hom-power} of $a \in A$ defined by \\
\par
$a^{1} = a$, \;\; $a^{n} = {\alpha}^{n-2}(a) a^{n-1}$ \hfill (3.10) \\
\\
for $n \geq 2$. In this setting of left Hom-powers, we have the following \\
\par
{\bf Theorem 3.10.} {\it Let $(L, \cdot , \alpha)$ be a Hom-Leibniz algebra and let $a \in L$. Then $L_{a^n} \circ \alpha = 0$, $n \geq 2$, where $L_{z}$ denotes the left multiplication by $z$ in $(L, \cdot , \alpha)$, i.e. $L_{z}x = z \cdot x$, $x \in L$}. \\
\par
{\bf Proof.} We proceed by induction on $n$ and the repeated use of Proposition 3.1(i).
\par
From Proposition 3.1(i), we get $a^{2}{\alpha}(z) = 0$, $\forall a, z \in L$ and thus the first step $n=2$ is verified. Now assume that, up to the degree $n$, we have $a^{n}{\alpha}(z) = 0$, $\forall a, z \in L$. Then Proposition 3.1(i) implies that \\ $ (a^{n}{\alpha}^{n-1}(a) + {\alpha}^{n-1}(a)a^{n}){\alpha}(z) =0$, i.e. $(a^{n}\alpha ({\alpha}^{n-2}(a)) + {\alpha}^{n-1}(a)a^{n}){\alpha}(z) =0$. The application of the induction hypothesis to $a^{n}\alpha ({\alpha}^{n-2}(a))$ leads to $({\alpha}^{n-1}(a)a^{n}){\alpha}(z) =0$, i.e. $({\alpha}^{(n+1)-2}(a)a^{(n+1)-1}){\alpha}(z) =0$ which means (by (3.10)) that $a^{n+1}{\alpha}(z) =0$. Therefore we conclude that $a^{n}{\alpha}(z) = 0$, $\forall n \geq 2$, i.e.  $L_{a^n} \circ \alpha = 0$, $n \geq 2$. \hfill $\square$ \\
\par
{\bf Remark 3.11.} We observed that Theorem 3.10 above is an $\alpha$-twisted version of a result of D.W. Barnes ([4], Theorem 1.2 and Corollary 1.3), related to left Leibniz algebras. Indeed, setting $\alpha = Id$, Theorem 3.10 reduces to the result of Barnes. \\
\\
{\bf References} \\
\\
$[1]$ M.A. Akivis. {\it Local algebras of a multidimensional three-web}, Siberian Math. J., {\bf 17} (1976), 3-8. \\
$[2]$ A.A. Albert. {\it On the power-associativity of rings}, Summa Brasil. Math., {\bf 2} (1948), 21-32. \\
$[3]$ Sh.A. Ayupov and B.A. Amirov. {\it On Leibniz algebras}, in: Algebras and Operator Theory, Proceedings of the Colloquium in Tashkent, Kluwer (1998), 1-13. \\
$[4]$ D.W. Barnes. {\it Engel subalgebras of Leibniz algebras}, arXiv:0810.2849v1. \\
$[5]$ J.T. Hartwig, D. Larsson and S.D. Silvestrov. {\it Deformations of Lie algebras using $\sigma$-derivations}, J. Algebra, {\bf 295} (2006), 314-361. \\
$[6]$ K. H. Hofmann and K. Strambach. {\it Lie's fundamental theorems for local analytical loops}, Pacific J. Math., {\bf 123} (1986), 301-327. \\
$[7]$ A.N. Issa. {\it Hom-Akivis algebras}, arXiv:1003.4770v3. \\
$[8]$ D. Larsson and S.D. Silvestrov. {\it Quasi-Lie algebras}, Contemp. Math., {\bf 391} (2005), 241-248. \\
$[9]$ J.-L. Loday. {\it Une version non commutative des alg\`ebres de Lie: les alg\`ebres de Leibniz}, Enseign. Math., {\bf 39} (1993), 269-293. \\
$[10]$ A. Makhlouf. {\it Hom-alternative algebras and Hom-Jordan algebras}, arXiv:0909.0326v1. \\
$[11]$ A. Makhlouf. {\it Paradigm of nonassociative Hom-algebras and Hom-superalgebras}, arXiv:1001.4240v1. \\
$[12]$ A. Makhlouf and S.D. Silvestrov. {\it Hom-algebra structures}, J. Gen. Lie Theory Appl., {\bf 2} (2008), 51-64. \\
$[13]$ D. Yau. {\it Enveloping algebras of Hom-Lie algebras}, J. Gen. Lie Theory Appl., {\bf 2} (2008), 95-108. \\
$[14]$ D. Yau. {\it Hom-algebras and homology}, J. Lie Theory, {\bf 19} (2009), 409-421. \\
$[15]$ D. Yau. {\it Hom-Novikov algebras}, arXiv:0909.0726v1. \\
$[16]$ D. Yau. {\it Hom-Maltsev, Hom-alternative and Hom-Jordan algebras}, arXiv:1002.3944v1. \\
$[17]$ D. Yau. {\it Hom-power associative algebras}, arXiv:1007.4118v1.
\end{document}